\documentclass{amsart}
\usepackage[utf8]{inputenc}
\usepackage{latexsym,amssymb,amsmath,amsthm,amscd,graphicx,enumerate}
\usepackage{amsfonts}
\usepackage{xcolor}
\usepackage{tikz}
\usepackage[english]{babel}
\usepackage[a4paper,bindingoffset=0.2in,%
            left=1in,right=1in,top=1in,bottom=1in,%
            footskip=.25in]{geometry}

\newtheorem{theorem}{Theorem}[section]
\newtheorem{lemma}[theorem]{Lemma}
\newtheorem{remark}[theorem]{Remark}
\newtheorem{cor}[theorem]{Corollary}
\newtheorem{defn}[theorem]{Definition}
\numberwithin{equation}{section}

\begin{document}

\title[Boundedness of classical operators on Generalized weighted Morrey spaces]{Boundedness of the Hardy-Littlewood
Maximal Operator, Fractional Integral Operators, and Calderon-Zygmund Operators on Generalized Weighted Morrey Spaces}

\author[Y.~Ramadana]{Yusuf Ramadana}
\address{Faculty of Mathematics and Natural Sciences, Bandung Institute of Technology, Bandung 40132, Indonesia}
\email{yusuframadana.96@gmail.com}

\author[H.~Gunawan]{Hendra Gunawan}
\address{Faculty of Mathematics and Natural Sciences, Bandung Institute of Technology, Bandung 40132, Indonesia}
\email{hgunawan@math.itb.ac.id}



\begin{abstract}
In this paper we investigate the boundedness of classical operators, namely the Hardy-Littlewood maximal operator,
fractional integral operators, and Calderon-Zygmund operators, on generalized weighted Morrey spaces
and generalized weighted weak Morrey spaces. We prove that the operators are bounded on these spaces, under certain
assumptions.
\end{abstract}

\maketitle

\section{Introduction}

For $a\in \mathbb{R}^n$ and $r>0$, we denote by $B(a,r)$ an open ball centered at $a$ with radius $r$.
For a set $E$ in $\mathbb{R}^n$, we denote by $E^{\rm c}$ the complement of $E$. Moreover, if $E$ is a
measurable set in $\mathbb{R}^n$, then $|E|$ denotes the Lebesgue measure of $E$.

The Hardy-Littlewood maximal operator $M$ and fractional maximal operator $M_\alpha$, where $0 \leq \alpha < n$,
are defined by
$$
Mf(x) := \sup_{r>0} \frac{1}{|B(x,r)|} \int_{B(x,r)} |f(y)| dy,\quad x\in\mathbb{R}^n
$$
and
$$
M_{\alpha}f(x) := \sup_{r>0} \frac{1}{|B(x,r)|^{1-\frac{\alpha}{n}}} \int_{B(x,r)} |f(y)| dy,\quad x\in\mathbb{R}^n
$$
for locally integrable functions $f$ on $\mathbb{R}^n$. It is a well-known property that $M$ is bounded on
$L^p=L^p(\mathbb{R}^n)$ for $1<p\le\infty$ and is bounded from $L^1$
to $L^{1,\infty}=L^{1,\infty}(\mathbb{R}^n)$, see e.g. \cite{Grafakos, Sawano}.

For $0< \alpha < n$, we also know the Riesz potential or the fractional integral operator $I_\alpha$ defined by
$$
I_{\alpha}f(x):= \int_{\mathbb{R}^n} \frac{f(y)}{|x-y|^{n-\alpha}} dy,\quad x\in \mathbb{R}^n
$$
for suitable functions $f$ on $\mathbb{R}^n$. The operator $I_\alpha$ is bounded from $L^p$ to $L^q$ for
$1<p<\infty$ and $\frac1q=\frac1p-\frac\alpha n$, see e.g. \cite{Stein93}. Since from the definitions we have
\begin{equation}\label{Eq1}
M_\alpha f(x) \leq C_n I_\alpha (|f|)(x),\quad x\in\mathbb{R}^n
\end{equation}
where $C_n$ is the Lebesgue measure of the unit ball in $\mathbb{R}^n,$ it thus follows that the operator $M_\alpha$
is also bounded from $L^p$ to $L^q$ for $1<p<\infty$ and $\frac1q=\frac1p-\frac\alpha n$.

Next, we consider the Calderon-Zygmund operator. Let $T = T_K$ be a linear operator from the Schwartz class $\mathcal{S}=
\mathcal{S}(\mathbb{R}^n)$ to $\mathcal{S}'$ which is $L^2-$bounded and for each $f \in \mathcal{C}_{\rm c}^{\infty}(\mathbb{R}^n)$ we have
$$
Tf(x) := \int_{\mathbb{R}^n} K(x,y) f(y) dy,\quad x \notin {\rm supp}(f)
$$
where $K=K(\cdot,\cdot)$ is the standard kernel defined on $\mathbb{R}^n \times \mathbb{R}^n$ except for the diagonal
$\{(x,x) : x\in \mathbb{R}^n \}$ with the following properties: there exists a constant $A>0$ for which
$$
|K(x, y)| \leq \frac{A}{|x-y|^n},\quad x \neq y,
$$
and for some $\delta > 0,$
$$
|K(x,y) - K(x',y)| \leq \frac{A |x-x'|^{\delta}}{(|x-y|+|x'-y|)^{n+\delta}},\quad |x-x'| \leq \frac{1}{2} \max (|x-y|, |x'-y|)
$$
and
$$
|K(x,y) - K(x,y')| \leq \frac{A |y-y'|^{\delta}}{(|x-y|+|x-y'|)^{n+\delta}},\quad |y-y'| \leq \frac{1}{2} \max (|x-y|, |x-y'|).
$$
The operator $T$ is called the Calderon-Zygmund operator, which was first introduced by Coifman and Meyer in 1979 \cite{Coifman}.
The operator is bounded on $L^p$ for $1<p<\infty$ and from $L^1$ to $L^{1,\infty}$ \cite{Grafakos}.

In this paper, we are interested in studying the boundedness of the above operators on generalized weighted Morrey spaces.
For $1\le p<\infty$ and $0\le\lambda<n$, the classical Morrey space $\mathcal{M}^{p, \lambda}= \mathcal{M}^{p, \lambda}(\mathbb{R}^n)$
was first introduced in \cite{Morrey} with the following norm
$$
\| f \|_{\mathcal{M}^{p, \lambda}} := \sup_{a \in \mathbb{R}^n, r>0} \frac{1}{r^\lambda} \left(\int_{B(a,r)} |f(x)|^p dx
\right)^{\frac{1}{p}} < \infty.
$$
The same space may be denoted by $\mathcal{M}^p_q = \mathcal{M}^p_q(\mathbb{R}^n)$ equipped with the norm
$$
\|f\|_{\mathcal{M}_p^q} := \sup_{a\in \mathbb{R}^n, r>0} \frac{1}{|B(a,r)|^{\frac{1}{p} - \frac{1}{q}}} \left(\int_{B(a,r) }
|f(x)|^p dx \right)^{\frac{1}{p}} =  \sup_{a\in \mathbb{R}^n, r>0} \frac{1}{|B(a,r)|^{\frac{1}{p} - \frac{1}{q}}}
\|f\|_{L^p (B(a,r))}
$$
where $1 \leq p \leq q < \infty$, as used widely in, for examples, \cite{Gunawan, Lin, Sawano}. Note that if we set $p=q$, then
$\mathcal{M}^p_q = L^p.$ In companion with $\mathcal{M}^p_q$, one may also define the weak Morrey space $W\mathcal{M}^p_q =
W\mathcal{M}^p_q(\mathbb{R}^n)$ as the set of all measurable functions $f$ on $\mathbb{R}^n$ such that
$$
\|f\|_{W\mathcal{M}^p_q} := \sup_{a \in \mathbb{R}^n, r>0} \frac{1}{|B(a,r)|^{\frac{1}{p} - \frac{1}{q}}} \|f\|_{WL^p(B(a,r))}<\infty
$$
where $\|f\|_{WL^p(B(a,r)}:=\sup\limits_{\gamma>0} \gamma\,|\{x\in B(a,r)\,:\,|f(x)|>\gamma\}|^{1/p}$.
The last two definitions were used in, for example, \cite{GunawanSawano}.

According to \cite{DiFazio}, $T$ is bounded on $\mathcal{M}_p^q$ for $1<p\leq q < \infty$ and is bounded from
$\mathcal{M}_1^q$ to $W \mathcal{M}_1^q$ for $1 \leq q < \infty$. In addition, $M$ is bounded on $\mathcal{M}_p^q$
for $1<p\leq q < \infty$ and is bounded from $\mathcal{M}_1^q$ to $W \mathcal{M}_1^q$ for $1 \leq q < \infty$ \cite{Chiarenza}.
Moreover, $I_\alpha$ is bounded from one Morrey space to another under certain conditions \cite{Adams, Peetre}.

In \cite{Mizuhara, Nakai}, the Morrey space $\mathcal{M}^p_q$ was generalized to $\mathcal{M}^p_\psi =
\mathcal{M}^p_\psi (\mathbb{R}^n)$, which consists of all locally integrable functions $f$ on $\mathbb{R}^n$ such that
the norm
$$
\|f\|_{\mathcal{M}^p_\psi} := \sup_{a \in \mathbb{R}^n} \left( \frac{1}{\psi (B(a,r))} \int_{\mathbb{R}^n} |f(x)|^p dx
\right)^{\frac{1}{p}} < \infty.
$$
Here $\psi$ is a function from $(0,\mathbb{R}^n \times \infty)$ to $(0,\infty)$ satisfying certain conditions.
Moreover, the weak generalized Morrey space $W \mathcal{M}^p_\psi = W \mathcal{M}^p_\psi (\mathbb{R}^n)$
where $1<p<\infty$ was defined as the set of all functions $f$ for which there exists a constant $C>0$ such that
$$
\frac{\gamma^p}{\psi(B)} |\{x \in B : |f(x)|>\gamma\}| \leq C
$$
for every ball $B=B(a,r)$ and $\gamma > 0.$ We can see that if we set $\psi(B(a,r)) = |B(a,r)|^{1- \frac{p}{q}}$
where $1 \leq q < \infty,$ then $\mathcal{M}^p_{\psi}=\mathcal{M}^p_q$. In \cite{Nakai}, Nakai investigated the
sufficient conditions on the function $\psi$ to ensure the boundedness of the operator $M, T$, and $I_\alpha$
on these spaces. The results are the following.

\begin{theorem}\label{Nakai1} \cite{Nakai}
Assume that there is a constant $C>0$ such that for any $a \in \mathbb{R}^n$ and $r>0$,
\begin{enumerate}[(a)]
	\item $r \leq t \leq 2r$ implies $\frac{1}{C} \leq \frac{\psi(a,t)}{\psi(a,r)} \leq C,$ and
	\item $\int_{r}^{\infty} \frac{\psi(a,t)}{t^{n+1}} dt \leq C \frac{\psi(a,r)}{r^n}.$
\end{enumerate}
Then,
\begin{enumerate}[(1)]
	\item $M$ and $T$ are bounded on $\mathcal{M}^p_\psi$ for $1<p<\infty$.
	\item $M$ and $T$ are bounded from $\mathcal{M}^1_\psi$ to $W \mathcal{M}^1_\psi.$
\end{enumerate}
\end{theorem}

\begin{theorem}\label{Nakai2} \cite{Nakai}
Let $0<\alpha<n, 1\leq p<n/\alpha, 1/q = 1/p-\alpha/n.$ Assume that for $r \leq t \leq 2r$ we have
$$
\frac{1}{C} \leq \frac{\psi(a,t)}{\psi(a,r)} \leq C,
$$
and also assume that for $r>0$ we have
$$
\int_r^\infty \frac{w(a,t)}{t^{n-\alpha p +1}} dt \leq C \frac{w(a,r)}{r^{n-\alpha p}}.
$$
\begin{enumerate}[(1)]
	\item If $p>1,$ then $I_\alpha$ is bounded from $\mathcal{M}^p_{\psi}$ to $\mathcal{M}^q_{\psi^{q/p}}.$
	\item If $p=1,$ then $I_\alpha$ is bounded from $\mathcal{M}^1_{\psi}$ to $W\mathcal{M}^q_{\psi^{1/q}}.$
\end{enumerate}
\end{theorem}

The results are sharper than those obtained by Mizuhara in \cite{Mizuhara}. Mizuhara used the assumption
that $\psi$ was a growth function satisfying doubling condition with a doubling constant $1 \leq D = D(\psi) < 2^n$.

In \cite{Guliyev}, Guliyev gave a similar definition for generalized Morrey spaces as Nakai and Mizuhara did.
Guliyev defined the norm on $\mathcal{M}^p_\phi = \mathcal{M}^p_\phi(\mathbb{R}^n)$ by
$$
\|f\|_{\mathcal{M}^p_\phi} := \sup_{a \in \mathbb{R}^n}\frac{r^{-\frac{n}{p}}}{\phi(a,r)} \|f\|_{L^p(B(a,r))}
$$
and the norm on $W \mathcal{M}^p_\phi = W \mathcal{M}^p_\phi (\mathbb{R}^n)$ by
$$
\|f\|_{W \mathcal{M}^p_\phi} = \sup_{a \in \mathbb{R}^n, r>0} \frac{r^{-\frac{n}{p}}}{\phi(a,r)} \|f\|_{WL^p(B(a,r))}.
$$

In the proof the Theorem \ref{Nakai1}, Nakai used the Hardy Littlewood Maximal operator $M$ and the assumption (a) played an
important role in the proof. In \cite{Guliyev}, the assumption (a) was removed. Guliyev also investigated the
boundedness of $M$ and $T$ between two generalized Morrey spaces $\mathcal{M}^p_{\phi_1}$ and $\mathcal{M}^p_{\phi_2}$
for $1<p<\infty,$ and from $\mathcal{M}^1_{\phi_1}$ to $W\mathcal{M}^1_{\phi_2}$ for some functions $\phi_1$ and $\phi_2$
on $\mathbb{R}^n \times (0,\infty)$ as stated in the following theorem.

\begin{theorem}\label{Guliyef} \cite{Guliyev}
Let $1 \leq p < \infty$ and the functions $\phi_1(a,r)$ and $\phi_2(a,r)$ satisfy
$$
\int_r^\infty \phi_1(a,t) \frac{dt}{t} \leq C \phi_2(a,r)
$$
for every $(a,r) \in \mathbb{R}^n \times (0,\infty)$ where $C$ does not depend on $a$ and $r.$
Then, $M$ and $T$ are bounded from $\mathcal{M}^p_{\phi_1}$ to $\mathcal{M}^p_{\phi_2}$ for $1<p<\infty$
and are bounded from $\mathcal{M}^1_{\phi_1}$ to $W\mathcal{M}^1_{\phi_2}$
\end{theorem}

Similar to Guliyev's definitions, we may also defined the generalized Morrey space $\mathcal{M}^p_\phi$
as the set of all locally integrable functions $f$ on $\mathbb{R}^n$ such that $\|f\|_{\mathcal{M}^p_\phi} < \infty$ where
$$
\begin{aligned}
\|f\|_{\mathcal{M}^p_{\phi}} &= \sup_{a \in \mathbb{R}^n, r>0} \frac{1}{\phi(a,r)} \left(\frac{1}{|B(a,r)|}
\int_{B(a,r)} |f(x)|^p dx\right)^{\frac{1}{p}} \\
& = \sup_{a \in \mathbb{R}^n, r>0} \frac{1}{\phi(a,r)} \cdot \frac{1}{|B(a,r)|^{\frac{1}{p}}} \|f\|_{L^p(B(a,r))}
\end{aligned}
$$
for $1 \leq p < \infty$ and a positive function $\phi$ on $\mathbb{R}^n \times (0,\infty).$
Moreover, we may also defined generalized weak Morrey space $\mathcal{M}^p_\phi$ as the set of all locally integrable functions
$f$ on $\mathbb{R}^n$ such that $\|f\|_{W\mathcal{M}^p_\phi} < \infty$ where
$$
\|f\|_{\mathcal{M}^p_\phi} = \sup_{a \in \mathbb{R}^n, r>0} \frac{1}{\phi(a,r)} \cdot \frac{1}{|B(a,r)|^{\frac{1}{p}}} \|f\|_{WL^p(B(a,r))}.
$$

The purpose of this article is to investigate the boundedness of $M, I_\alpha,$ and $T$ on generalized weighted Morrey spaces 
and weighted generalized weak Morrey spaces. 
The definition of generalized weighted Morrey spaces and generalized weighted weak Morrey spaces is formulated by replacing 
$|B(a,r)|$ with $w(B(a,r)):=\int_{B(a,r)} w(x)dx$ for some weight function $w$ which we shall discuss in the next section. 
The results we obtain thus generalize the previous results.

\section{$A_p$ Weights}

In this section, we discuss the $A_p$ weights and the weighted Lebesgue spaces. We also present the
definition of the generalized weighted Morrey spaces, the generalized weighted weak Morrey spaces, and some lemmas that we shall
use to prove the main results about the boundedness of the three classical operators on the generalized weighted Morrey spaces and
the generalized weighted weak Morrey spaces.

A weight $w$ is a nonnegative locally integrable function on $\mathbb{R}^n$ taking values in the interval $(0,\infty)$ almost
everywhere. The weight class that we use in this article is the Muckenhoupt class $A_p.$

\begin{defn} \cite{Garcia}
For $1 < p < \infty,$ we denote by $A_p$ the set of all weights $w$ on $\mathbb{R}^n$ for which there exists a constant $C>0$
such that
$$
\left(\frac{1}{|B(a,r)|} \int_{B(a,r)} w(x) dx\right) \left(\frac{1}{|B(a,r)|} \int_{B(a,r)} w(x)^{-\frac{1}{p-1}} dx\right)^{p-1}
\leq C
$$
for every ball $B(a,r)$ in $\mathbb{R}^n$. For $p=1$, we denote by $A_1$ the set of all weights $w$ for which there exists a
constant $C>0$ such that
$$
\frac{1}{|B(a,r)|} \int_{B(a,r)} w(x) dx \leq C \|w\|_{L^{\infty}(B(a,r))}
$$
for every ball $B(a,r)$ in $\mathbb{R}^n$.
\end{defn}

\begin{remark}
The last inequality is equivalent to the following
$$
\left(\frac{1}{|B(a,r)|} \int_{B(a,r)} w(x) dx)\right) \cdot \|w^{-1}\|_{L^{\infty}(B(a,r))} \leq C
$$
for every ball $B(a,r)$ in $\mathbb{R}^n$.
\end{remark}

\begin{theorem}\label{garcia}\cite{Garcia}
For each $1\leq p < \infty$ and $w \in A_p$, there exists $C>0$ such that
$$
\frac{w(B)}{w(E)} \leq C \left(\frac{|B|}{|E|}\right)^p
$$
for every ball $B$ and measurable sets $E \subseteq B$ where
$$
w(B) = \int_B w(x) dx.
$$
\end{theorem}

Associated to a weight function $w\in A_p$ with $1 \leq p < \infty$ and a measurable set $\Omega$ in $\mathbb{R}^n$, we define the weighted Lebesgue space $L^{p,w}(\Omega)$
to be the set of all measurable functions $f$ on $\Omega$ for which
$$
\|f\|_{L^{p,w}(\Omega)} := \left(\int_{\Omega} |f(x)|^p w(x) dx \right)^\frac1p < \infty.
$$
In addition, we denote by $WL^{p,w}(\Omega)$ the weighted weak Lebesgue space that consists of all measurable functions $f$
on $\Omega$ for which
$$
\|f\|_{WL^{p,w}(\Omega)} := \sup_{\gamma>0} \gamma w(\{x\in \Omega : |f(x)|>\gamma\})^{\frac{1}{p}} < \infty.
$$
We write $L^{p,w} = L^{p,w}(\mathbb{R}^n)$ and $WL^{p,w} = WL^{p,w}(\mathbb{R}^n).$ Notice that if $w$ is a constant function a.e., then we have that $L^{p,w} = L^p$ and $WL^p = L^{p,\infty}$.
We note from \cite{Garcia} that $w \in A_p$ if and only if $M$ is bounded on $L^{p,w}$ for $1<p<\infty$ and $w \in A_1$ if
and only if $M$ is bounded from $L^{1,w}$ to $WL^{1,w}$.

Related to the fractional integral operator $I_\alpha,$ we have another kind of class of weights $A_{p,q}$.

\begin{defn} \cite{Wheeden, Sawano}
Let $1 < p < q < \infty$ and $p'$ satisfies $1/p+1/p'=1.$ We denote by $A_{p,q}$ the collection of all weight functions $w$
satisfying
$$
\left(\frac{1}{|B(a,r)|} \int_{B(a,r)} w(x)^q dx\right)^{\frac{1}{q}} \left(\frac{1}{|B(a,r)|}
\int_{B(a,r)} w(x)^{-p'}\right)^{\frac{1}{p'}} \leq C
$$
for every $(a,r)\in \mathbb{R}^n \times (0,\infty)$ where $C$ is a constant independent of $a$ and $r.$
For $p=1$ and $q>1$, we denote by $A_{1,q}$ the collection of weight functions $w$ for which there exists a constant $C>0$ such that
$$
\left(\frac{1}{|B(a,r)|} \int_{B(a,r)} w(x)^q dx\right)^{1/q} \leq C \|w\|_{L^\infty(B(a,r))}
$$
for every $(a,r)\in \mathbb{R}^n\times (0,\infty)$.
\end{defn}

\begin{theorem}\label{thm2.5}
Let $1 \leq p < q.$ Then $w \in A_{p,q}$ if and only if $w^q \in A_{q/p'+1}$. Moreover, if $w\in A_{p,q}$, then $w^p\in A_p$ and $w^q \in A_q$.
\end{theorem}

\textit{Proof.} Let $p>1.$ Observe that $w\in A_{p,q}$ means that
$$
\left(\frac{1}{|B(a,r)|} \int_{B(a,r)} w(y)^q dy\right)^\frac1q \left(\frac{1}{|B(a,r)|}
\int_{B(a,r)}w(y)^{-p'} dy \right)^{-\frac{1}{p'}} \leq C
$$
for every $(a,r) \in \mathbb{R}^n\times (0,\infty).$ Raising both sides to the $q-$th power, we obtain
$$
\left(\frac{1}{|B(a,r)|} \int_{B(a,r)} w(y)^q dy\right) \left(\frac{1}{|B(a,r)|} \int_{B(a,r)}w(y)^{-q'} dy \right)^{-\frac{q}{q'}} \leq C
$$
for every $(a,r) \in \mathbb{R}^n\times (0,\infty).$ This means that $w^q\in A_{q/p'+1}$.

For $p=1$ and $q>1$, we have $p'=\infty$ and $\frac{q}{p'}+1=1.$ Here $w\in A_{1,q}$ means that
$$
\frac{1}{|B(a,r)|} \int_{B(a,r)} w(x)^q dx \leq C \|w^q\|_{L^\infty (B(a,r))}
$$
for every $(a,r) \in \mathbb{R}^n\times (0,\infty).$ This tells us precisely that $w^q\in A_1.$

For the second part, first suppose that $w\in A_{p,q},$ where $1<p<q$. Then,
$$
\left(\frac{1}{|B(a,r)|} \int_{B(a,r)} w(y)^q dy\right)^\frac1q \left(\frac{1}{|B(a,r)|}
\int_{B(a,r)}w(y)^{-p'} dy \right)^{-\frac{1}{p'}} \leq C
$$
for every $(a,r)\in \mathbb{R}^n\times (0,\infty).$ By H\"older's inequality,
$$
\left(\frac{1}{|B(a,r)|} \int_{B(a,r)} w(y)^p dy\right)^\frac1p \leq \left(\frac{1}{|B(a,r)|}
\int_{B(a,r)} w(y)^q dy\right)^\frac1q
$$
which implies that
$$
\left(\frac{1}{|B(a,r)|} \int_{B(a,r)} w(y)^p dy\right)^\frac1p \left(\frac{1}{|B(a,r)|}
\int_{B(a,r)}w(y)^{-p'} dy \right)^{-\frac{1}{p'}} \leq C
$$
or
$$
\left(\frac{1}{|B(a,r)|} \int_{B(a,r)} w(y)^p dy\right) \left(\frac{1}{|B(a,r)|} \int_{B(a,r)}w(y)^{-\frac{p}{p-1}} dy
\right)^{\frac{1}{p-1}} \leq C
$$
for every $(a,r)\in\mathbb{R}^n\times (0,\infty),$ which tells us that $w^p\in A_p.$

Next, if $w\in A_{1,q}$ for $1<q$, then
$$
\left(\frac{1}{|B(a,r)|} \int_{B(a,r)} w(x)^q dx\right)^\frac1q \leq C \|w\|_{L^\infty (B(a,r))}
$$
for every $(a,r)\in \mathbb{R}^n\times (0,\infty).$ By H\"older's inequality, we have
$$
\frac{1}{|B(a,r)|} \int_{B(a,r)} w(x) dx \leq C \|w\|_{L^\infty (B(a,r))}.
$$
This means that $w\in A_1.$

Finally, to show that $w^q\in A_q$ for $1\le p<q$, we observe that $q'<p'$, and hence
$$
\begin{aligned}
	& \left(\frac{1}{|B(a,r)|} \int_{B(a,r)} w(y)^q dy\right)^\frac1q \left(\frac{1}{|B(a,r)|}
\int_{B(a,r)}w(y)^{-q'} dy \right)^{\frac{1}{q'}} \\
	& \qquad\leq \left(\frac{1}{|B(a,r)|} \int_{B(a,r)} w(y)^q dy\right)^\frac1q \left(\frac{1}{|B(a,r)|}
\int_{B(a,r)}w(y)^{-p'} dy \right)^{\frac{1}{p'}} \leq C.
\end{aligned}
$$
One should interpret the above inequality appropriately for the case where $p=1$. Therefore, we obtain
$$
\left(\frac{1}{|B(a,r)|} \int_{B(a,r)} w(y)^q dy\right)\left(\frac{1}{|B(a,r)|}
\int_{B(a,r)}(w(y)^q)^{-\frac{1}{q-1}} dy \right)^{q-1} \leq C,
$$
whence $w^q\in A_q$, and so the proof is complete.
\qed

On weighted Lebesgue spaces, we rewrite the following results of Garcia-Cuerva and Rubio de Francia \cite{Garcia} (for $M$ and $T$),
Muckenhoupt and Wheeden \cite{Wheeden} (for $I_\alpha$), and Quek \cite{Quek} and Yabuta \cite{Yabuta} (for $T$):

\begin{theorem}\label{thm3.1} \cite{Garcia}
Let $1<p<\infty.$ Then, $M$ is bounded on $L^{p,w}$ if $w\in A_p.$ Moreover, $M$ is bounded from
$L^{1,w}$ to $WL^{1,w}$ if $w\in A_1$.
\end{theorem}

\begin{theorem}\label{thm4.1} \cite{Wheeden}
Let $0<\alpha<n, 1<p<n/\alpha, 1/q=1/p-\alpha/n,$ and $w\in A_{p,q}$. Then, the operator $I_\alpha$ is bounded from $L^{p,w^p}$
to $L^{q,w^q}$. Moreover, if $w\in A_{1,p}$ with $1/q = 1-\alpha/n,$ then $I_\alpha$ is bounded from $L^{1,w}$ to $WL^{q,w^q}$.
\end{theorem}

\begin{theorem}\label{thm5.1} \cite{Garcia, Quek, Yabuta}
Let $1\leq p< \infty.$ Then, $T$ is bounded from $L^{p,w}$ if $w \in A_p$ and bounded from
$L^{1,w}$ to $WL^{1,w}$ if $w\in A_1$.
\end{theorem}

For $1 \leq p < \infty$ and a weight function $w$ on $\mathbb{R}^n$, the space $L^{p,w}_{loc}$ collects all measurable functions $f$ such that $\mathcal{X}_B \cdot f$ is in $L^{p,w}$ for each ball $B$ in $\mathbb{R}^n$. We now present the definition of the generalized weighted Morrey spaces and the generalized weighted weak Morrey spaces
which will become the spaces of our interest in this article.

\begin{defn}\label{wgms}
Let $1 \leq p <\infty, w \in A_p,$ and $\psi$ be a nonnegative function on $\mathbb{R}^n \times (0,\infty)$.
The {\it generalized weighted Morrey space} $\mathcal{M}^{p,w}_\psi = \mathcal{M}^{p,w}_\psi(\mathbb{R}^n)$ is the set of
all functions $f\in L^{p,w}_{loc}$ such that
$$
\begin{aligned}
\|f\|_{\mathcal{M}^{p,w}_\psi}
&= \sup_{a\in \mathbb{R}^n, r>0} \frac{1}{\psi (a,r)} \left(\frac{1}{w(B(a,r))}
\int_{B(a,r)} |f(x)|^p w(x) dx\right)^{1/p} \\
&=\sup_{a\in \mathbb{R}^n, r>0} \frac{1}{\psi(a,r)} \frac{1}{w(B(a,r))^{\frac{1}{p}}} \|f\|_{L^{p,w}(B(a,r))}<\infty.
\end{aligned}
$$
\end{defn}

\begin{defn}\label{wgwms}
Let $1 \leq p < \infty, w \in A_p,$ and $\psi$ be a nonnegative function on $\mathbb{R}^n \times (0,\infty).$
The {\it generalized weighted weak Morrey space} $W \mathcal{M}^{p, w}_\psi = W \mathcal{M}^{p, w}_\psi (\mathbb{R}^n)$
is the set of all functions $f\in L^{p,w}_{loc}$ such that
$$
\begin{aligned}
\|f\|_{W \mathcal{M}^{p, w}_\psi}
&= \sup_{a \in \mathbb{R}^n, r>0} \sup_{\gamma>0} \frac{1}{\psi (a,r)} \frac{\gamma}{w(B(a,r))^{\frac{1}{p}}}
w(\{x\in\mathbb{R}^n : |f(x)|>\gamma\})^{\frac{1}{p}} \\
&= \sup_{a \in \mathbb{R}^n, r>0} \frac{1}{\psi (a,r)} \frac{1}{w(B(a,r))^{\frac{1}{p}}} \|f\|_{WL^{p,w}(B(a,r))}<\infty.
\end{aligned}
$$
\end{defn}

Note that if we set $w$ to be constant a.e., then $\mathcal{M}^{p,w}_\psi = \mathcal{M}^{p}_\psi$ and $W\mathcal{M}^{p,w}_\psi =
W\mathcal{M}^{p}_\psi$. Moreover, if we set $\psi(a,r) = |B(a,r)|^{-\frac{1}{q}}$ and $w$ is constant a.e., then
$\mathcal{M}^{p,w}_\psi$ is the classical Morrey space and $W\mathcal{M}^{p,w}_\psi$ is the classical weak Morrey space.

With the definitions \ref{wgms} and \ref{wgwms}, we shall investigate the boundedness of the classical operators:
the Hardy-Littlewood maximal operator, the fractional integral operators, the fractional maximal operators, and the Calderon-Zygmund
operators on those spaces in the next section.

We end this section with lemmas which will be used later in proving our main theorems.

\begin{lemma}\label{lem2.8}
Let $\varphi$ be a nonnegative function on $\mathbb{R}^n \times (0,\infty)$ such that the map $r \mapsto \varphi(a,r)$
is increasing for each $a\in \mathbb{R}^n.$ Let $w\in A_p$ where $1 \leq p < \infty$. Then, for every ball $B(a,r)$, we have
$$
\varphi(a,r) \leq C w(B(a,r))^{\frac{1}{p}} \sup_{r < s < \infty} \frac{1}{w(B(a,s))^{\frac{1}{p}}} \varphi(a,s)
$$
and
$$
\varphi(a,r) \leq C w(B(a,r))^{\frac{1}{p}} \int_r^\infty \frac{1}{w(B(a,s))^{\frac{1}{p}}} \varphi(a,s) \frac{ds}{s}
$$
where $C>0$ is independent of $a\in \mathbb{R}^n$, $r>0,$ and $\varphi$.
\end{lemma}

\textit{Proof.} Let $a\in \mathbb{R}^n$ and $r>0$. By Theorem \ref{garcia} and the fact that the $ s \mapsto \varphi (a,s)$ map
is increasing for each $a\in\mathbb{R}^n$, we have
$$
\begin{aligned}
w(B(a,r))^{\frac{1}{p}}  \sup_{r < s < \infty} \frac{1}{w(B(a,s))^{\frac{1}{p}}} \varphi(a,s)
& = \sup_{r < s < \infty} \frac{w(B(a,r))^{\frac{1}{p}}}{w(B(a,s))^{\frac{1}{p}}} \varphi(a,s) \\
& \geq \sup_{r < s < \infty} C \frac{|B(a,r)|}{|B(a,s)|} \varphi (a,s) \\ &= C \sup_{r < s < \infty} \frac{r^n}{s^n} \varphi(a,s) \\
& \geq C \sup_{r < s < \infty} \frac{r^n}{s^n} \varphi(a,r) \\
& = C \varphi(a,r).
\end{aligned}
$$
Moreover,
$$
\begin{aligned}
w(B(a,r))^{\frac{1}{p}} \int_r^\infty \frac{1}{w(B(a,s))^{\frac{1}{p}}} \varphi(a,s) \frac{ds}{s}
& = \int_r^\infty \frac{w(B(a,r))^{\frac{1}{p}}}{w(B(a,s))^{\frac{1}{p}}} \varphi (a,s) \frac{ds}{s} \\
& \geq \int_r^\infty C \frac{|B(a,r)|}{|B(a,s)|} \varphi(a,s) \frac{ds}{s} \\
&= C \int_r^\infty \frac{r^n}{s^n} \varphi(a,s) \frac{ds}{s} \\ & \geq C \int_r^{2r} \frac{r^n}{s^n} \varphi(a,r) \frac{ds}{s} \\
& \geq C \int_r^{2r} \frac{r^n}{(2r)^n} \varphi(a,r)\frac{ds}{2r} \\ &= C \varphi(a,r).
\end{aligned}
$$
Therefore,
$$
\varphi(a,r) \leq C w(B(a,r))^{\frac{1}{p}} \sup_{r < s < \infty} \frac{1}{w(B(a,s))^{\frac{1}{p}}} \varphi (a,s)
$$
and
$$
\varphi(a,r) \leq C w(B(a,r))^{\frac{1}{p}} \int_r^\infty \frac{1}{w(B(a,s))^{\frac{1}{p}}} \varphi (a,s) \frac{ds}{s},
$$
which proves the lemma.\qed

\begin{lemma}\label{lem2.9}
Let $1\leq p<\infty$ and $w\in A_p.$ For each $(a,s)\in\mathbb{R}^n \times (0,\infty)$ and $f\in L^{p,w}_{loc}$, we have
$$
\frac{1}{|B(a,s)|} \int_{B(a,s)} |f(y)| dy \leq \frac{C}{w(B(a,s))^\frac{1}{p}} \|f\|_{L^{p,w}(B(a,s))}
$$
where $C$ is constant independent of $a$, $s,$ and $f.$
\end{lemma}

\textit{Proof.} First, we assume that $1<p<\infty.$ Thus, by using H\"older's inequality and assumption that $w\in A_p,$ we have
for every $a\in \mathbb{R}^n$ and $s>0$,
$$
\begin{aligned}
& \frac{1}{|B(a,s)|} \int_{B(a,s)} |f(y)| dy \\
& = \frac{w(B(a,s))^\frac{1}{p}}{|B(a,s)|} \frac{1}{w(B(a,s))^{\frac{1}{p}}} \int_{B(a,s)} |f(y)|
\frac{w(y)^{\frac{1}{p}}}{w(y)^{\frac{1}{p}}} dy \\
& \leq \frac{w(B(a,s))^{\frac{1}{p}}}{|B(a,s)|} \frac{1}{w(B(a,s))^{\frac{1}{p}}} \left(\int_{B(a,s)}
|f(y)|^p w(y) dy\right)^{\frac{1}{p}}\left(\int_{B(a,s)} w(y)^{-\frac{p'}{p}} dy\right)^{\frac{1}{p'}}\\
& = \left(\frac{1}{|B(a,s)|} \int_{B(a,s)} w(y) dy\right)^{\frac{1}{p}} \left(\frac{1}{|B(a,s)|} \int_{B(a,s)}
w(y)^{-\frac{p'}{p}} dy \right)^{\frac{1}{p'}} \frac{1}{w(B(a,s))^{\frac{1}{p}}} \left(\int_{B(a,s)} |f(y)|^p w(y) dy\right)^{\frac{1}{p}} \\
& \leq \frac{C}{w(B(a,s))^{\frac{1}{p}}} \left(\int_{B(a,s)} |f(y)|^p w(y) dy\right)^{\frac{1}{p}} \\
& = \frac{C}{w(B(a,s))^\frac{1}{p}} \|f\|_{L^{p,w}(B(a,s))}.
\end{aligned}
$$
For $p=1,$ by using H\"older's inequality and assumption that $w\in A_1,$ we get
$$
\begin{aligned}
& \frac{1}{|B(a,s)|} \int_{B(a,s)} |f(y)| dy \\
& = \frac{w(B(a,s))}{|B(a,s)|} \frac{1}{w(B(a,s))} \int_{B(a,s)} |f(y)| \frac{w(y)}{w(y)} dy \\
& \leq \frac{w(B(a,s))}{|B(a,s)|} \frac{1}{w(B(a,s))} \int_{B(a,s)} |f(y)| w(y) dy \|w^{-1}\|_{L^\infty (B(a,s))} \\
& = \frac{1}{|B(a,s)|} \int_{B(a,s)} w(y) dy \|w^{-1}\|_{L^\infty (B(a,s))} \frac{1}{w(B(a,s))} \int_{B(a,s)} |f(y)| w(y) dy \\
& \leq \frac{C}{w(B(a,s))} \int_{B(a,s)} |f(y)| w(y) dy \\
& = \frac{C}{w(B(a,s))} \|f\|_{L^{1,w}(B(a,s))},
\end{aligned}
$$
as desired.\qed

\begin{cor}
For each $1 \leq p < \infty$ and $w\in A_p$, there exists $C>0$ such that for every $a\in \mathbb{R}^n$ and $r>0,$
we have
$$
\int_{\mathbb{R}^n \setminus B(a,r)} \frac{|f(y)|}{|a-y|^n} dy \leq C \int_r^\infty w(B(a,s))^{-\frac{1}{p}}
\|f\|_{L^{p,w}(B(a,s))} \frac{ds}{s}, \quad f \in L^{p,w}_{loc}.
$$
\end{cor}

\textit{Proof.} Let $a\in \mathbb{R}^n$ and $r>0$. By the Fubini's Theorem,
$$
\begin{aligned}
\int_{\mathbb{R}^n\setminus B(a,r)} \frac{f(y)}{|a-y|^n} dy
& = \int_{B(a,r)^c} |f(y)|
\int_{|a-y|}^\infty \frac{1}{s^n} \frac{ds}{s} dy \\
& = \int_r^\infty \int_{B(a,s)\setminus B(a,r)} |f(y)| dy \frac{1}{s^n} \frac{ds}{s} \\
& = C \int_r^\infty \frac{1}{|B(a,s)|} \int_{B(a,s)\setminus B(a,r)} |f(y)| dy \frac{ds}{s} \\
& \leq C  \int_r^\infty \frac{1}{|B(a,s)|} \int_{B(a,s)} |f(y)| dy \frac{ds}{s}.		
\end{aligned}
$$
Hence, by Lemma \ref{lem2.9},
$$
\begin{aligned}
\int_{\mathbb{R}^n \setminus B(a,r)} \frac{|f(y)|}{|a-y|^n} dy
& \leq C \int_r^\infty \frac{1}{|B(a,s)|} \int_{B(a,s)} |f(y)| dy \frac{ds}{s} \\
& \leq C \int_r^\infty \frac{1}{w(B(a,s))^{\frac{1}{p}}} \|f\|_{L^{p,w}(B(a,s))} \frac{ds}{s},
\end{aligned}$$
as claimed.
\qed

\section{Hardy-Littlewood Maximal Operator on Generalized Weighted Morrey Spaces}

In this section, we prove the boundedness of the Hardy-Littlewood operator $M$ on generalized weighted Morrey spaces and
generalized weighted weak Morrey spaces. Keeping in mind Theorem \ref{thm3.1} \cite{Garcia}, we have the following result.

\begin{theorem}\label{thm3.2}
Let $1\leq p < \infty, w \in A_p$. Then, for every $a \in \mathbb{R}^n$ and $r>0,$
$$
\|Mf\|_{L^{p,w}(B(a,r))} \leq C_1 w(B(a,r))^{\frac{1}{p}} \sup_{r < t < \infty} w(B(a,t))^{-\frac{1}{p}} \|f\|_{L^{p,w}B(a,t)},\quad
f \in L_{loc}^{p,w}
$$
for $1<p<\infty,$ and
$$
\|Mf\|_{WL^{1,w}(B(a,r))} < C_2 w(B(a,r)) \sup_{r < t < \infty} w(B(a,t))^{-1} \|f\|_{L^{1,w}B(a,t)},\quad f\in L_{loc}^{1,w}
$$
where $C_1$ and $C_2$ are constants that do not depend on $f, a,$ and $r.$
\end{theorem}

\textit{Proof.} Let $a\in \mathbb{R}^n$ and $r>0$, and write $f$ in the form of $f := f_1 + f_2$ where $f_1 :=
f\cdot \mathcal{X}_{B(a,2r)}.$ Assume that $1<p<\infty.$ Since $w\in A_p, M$ is bounded on $L^{p,w}.$ Thus,
$$
\|Mf\|_{L^{p,w} (B(a,r))} \leq \|Mf_1\|_{L^{p,w}(B((a,r))} + \|Mf_2\|_{L^{p,w}(B((a,r))},
$$
and
$$
\|Mf_1\|_{L^{p,w}(B(a,r))} \leq \|Mf_1\|_{L^{p,w}} \leq C \|f_1\|_{L^{p,w}} \leq C \|f\|_{L^{p,w}(B(a,2r))}.
$$
We can see that the map $r\mapsto \|f\|_{L^{p,w}(B(a,2r))}$ is increasing for each $a\in\mathbb{R}^n$. Then, by
Theorem \ref{garcia} and Lemma \ref{lem2.8},
$$
\|Mf_1\|_{L^{p,w} (B(a,r))} \leq C w(B(a,r))^{\frac{1}{p}} \sup_{r < t < \infty} \frac{1}{w(B(a,t))^{\frac{1}{p}}} \|f\|_{L^{p,w}(B(a,t))}.
$$

Let $x\in B(a,r).$ If $y\in B(x,t) \cap B(a,2r)^c$, then
$$
r=2r-r\leq |y-a| - |a-x| \leq |y-x| < t.
$$
In other words,
$$
\int_{B(x,t) \cap B(a,2r)^c} |f(y)| dy =0, \quad t \leq r.
$$
Moreover,
$$
|y-a| \leq |y-x|+|x-a| \leq t + r < 2t.
$$
It then follows that
$$
\begin{aligned}
Mf_2(x) & = \sup_{t>0} \frac{1}{|B(x,t)|} \int_{B(x,t)} |f_2 (y)| dy \\
& \leq \max \left(\sup_{t>r} \frac{1}{|B(x,t)|} \int_{B(x,t) \cap B(a,2r)^c} |f(y)| dy,
\sup_{0<t\leq r} \frac{1}{|B(x,t)|} \int_{B(x,t) \cap B(a,2r)^c} |f(y)| dy\right) \\
& =\sup_{t>r} \frac{1}{|B(x,t)|} \int_{B(x,t) \cap B(a,2r)^c} |f(y)| dy\\
& \leq \sup_{t> r} \frac{1}{|B(x,t)|} \int_{B(a,2t)} |f(y)| dy \\
& = C \sup_{t> 2r} \frac{1}{|B(a,t)|} \int_{B(a,t)} |f(y)| dy \\
& \leq C \sup_{t> r} \frac{1}{w(B(a,t))^{\frac{1}{p}}} \|f\|_{L^{p,w}(B(a,t))}.
\end{aligned}
$$
Hence,
$$
\|Mf_2\|_{L^{p,w} (B(a,r))} \leq C w(B(a,r))^{\frac{1}{p}} \sup_{r < t < \infty}
w(B(a,t))^{-\frac{1}{p}} \|f\|_{L^{p,w}B(a,t)},
$$
and so we conclude that
$$
\|Mf\|_{L^{p,w}(B(a,r))}\leq C_1 w(B(a,r))^{\frac{1}{p}} \sup_{r < t < \infty} w(B(a,t))^{-\frac{1}{p}} \|f\|_{L^{p,w} B(a,t)}.
$$
Assume now that $p=1.$ Thus,
$$
\|Mf\|_{WL^{1,w}(B(a,r))} \leq \|Mf_1\|_{L^{1,w}(B((a,r))} + \|Mf_1\|_{L^{1,w}(B((a,r))},
$$
and
$$
\|Mf_1\|_{WL^{1,w}(B(a,r))} \leq \|Mf_1\|_{L^{1,w}} \leq C \|f_1\|_{L^{1,w}} \leq C \|f\|_{L^{1,w}(B(a,2r))}.
$$
Since $M$ is bounded from $L^{1,w}$ to $WL^{1,w},$ we have
$$
\|Mf_1\|_{WL^{1,w}(B(a,r))} \leq \|Mf_1\|_{WL^{1,w}}\leq C \|f_1\|_{L^{1,w}} \leq C \|f\|_{L^{1,w}(B(a,2r))}.
$$
By Theorem \ref{garcia} and Lemma \ref{lem2.8},
$$
\|Mf_1\|_{WL^{1,w}(B(a,r))} \leq C w(B(a,r)) \sup_{r < t < \infty} w(B(a,t))^{-1} \|f\|_{L^{1,w}B(a,t)}.
$$
Since
$$
Mf_2(x) \leq C \sup_{r < t <\infty} w(B(a,t))^{-1} \|f\|_{L^{p,w}(B(a,t))}, \quad x\in B(a,r),
$$
by H\"older inequality and assumption that $w\in A_1,$ we obtain
$$
\begin{aligned}
& \|Mf_2\|_{WL^{1,w}(B(a,r))} \\ & = \sup_{\gamma>0} \gamma w(\{x\in B(a,r):|Mf_2 (x)| > \gamma\}) \\
& = \sup_{\gamma>0} \gamma \int_{\{x\in B(a,r):|Mf_2 (x)|>\gamma\}} |Mf_2 (x)|w(x)dx \cdot
\left\|\frac{1}{Mf_2}\right\|_{L^\infty (\{x\in B(a,r):|Mf_2(x)|>\gamma\})} \\
& \leq \sup_{\gamma>0} \gamma \int_{\{x \in B(a,r):|Mf_2(x)|>\gamma\}} |Mf_2(x)| w(x) dx \cdot \frac{1}{\gamma} \\
& \leq \int_{B(a,r)} |Mf_2(x)|w(x) dx \\
& \leq C \int_{B(a,r)} \sup_{r < t <\infty} w(B(a,t))^{-1} \|f\|_{L^{p,w}(B(a,t))} w(x) dx \\
& = C w(B(a,r)) \sup_{r < t <\infty} w(B(a,t))^{-1} \|f\|_{L^{p,w}(B(a,t))}.
\end{aligned}
$$
Therefore,
$$
\|Mf\|_{WL^{1,w} (B(a,r))} \leq C_2 \sup_{r < t <\infty} w(B(a,t))^{-1} \|f\|_{L^{p,w}(B(a,t))}
$$
and this proves the Theorem \ref{thm3.2}. \qed

\begin{theorem}\label{thm3.3}
Let $1 \leq p < \infty, w \in A_p$, and $M$ be the Hardy-Littlewood maximal operator. Suppose that $\psi_1$ and $\psi_2$
are two nonnegative functions on $\mathbb{R}^n \times (0,\infty)$ satisfying
$$
\sup_{r < t < \infty} \psi_1(a,t) \leq C \psi_2 (a,r)
$$
for every $(a,r)\in \mathbb{R}^n \times (0, \infty)$, where $C$ is a constant that does not depend on $a$ and $r.$ Then,
\begin{enumerate}[(1)]
	\item $M$ is bounded from $\mathcal{M}^{p,w}_{\psi_1}$ to $\mathcal{M}^{p,w}_{\psi_2}$ for $1<p<\infty.$
	\item $M$ is bounded from $\mathcal{M}^{1,w}_{\psi_1}$ to $W\mathcal{M}^{1,w}_{\psi_2}.$
\end{enumerate}
\end{theorem}

\textit{Proof.} First, assuming that $1<p<\infty,$ let $f\in \mathcal{M}^{p,w}_{\psi_1}.$ By using Theorem \ref{thm3.2}
and the hypothesis about $\psi_1$ and $\psi_2,$ we get
$$
\begin{aligned}
\|Mf\|_{\mathcal{M}^{p,w}_{\psi_2}}
&= \sup_{a\in \mathbb{R}^n, r>0} \frac{1}{\psi_2(a,r)} \left(\frac{1}{w(B(a,r))} \int_{B(a,r)} |Mf(x)|^p w(x) dx\right)^{\frac{1}{p}} \\
&= \sup_{a\in \mathbb{R}^n, r>0} \frac{1}{\psi_2(a,r)} w(B(a,r))^{-\frac{1}{p}} \|Mf\|_{L^{p,w}(B(a,r))} \\
&\leq C \sup_{a\in \mathbb{R}^n, r>0} \frac{1}{\psi_2(a,r)} \sup_{r < t < \infty} w(B(a,t))^{-\frac{1}{p}} \|f\|_{L^{p,w} B(a,t)}\\
&= C \sup_{a\in \mathbb{R}^n, r>0} \frac{1}{\psi_2(a,r)} \sup_{r < t <\infty} \frac{\psi_1(a,t)}{\psi_1(a,t)}
w(B(a,t))^{-\frac{1}{p}} \|f\|_{L^{p,w}(B(a,t))} \\
&\leq C \sup_{a\in \mathbb{R}^n, r>0} \frac{1}{\psi_2(a,r)} \sup_{r < t < \infty} \psi_1 (a,t) \|f\|_{\mathcal{M}^{p,w}_{\psi_1}} \\
&= C \|f\|_{\mathcal{M}^{p,w}_{\psi_1}} \sup_{a \in \mathbb{R}^n, r>0} \frac{1}{\psi_2(a,r)} \sup_{r < t < \infty} \psi_1(a,t) \\
&\leq C \|f\|_{\mathcal{M}^{p,w}_{\psi_1}}.
\end{aligned}
$$
Therefore, we conclude that $M$ is bounded from $\mathcal{M}^{p,w}_{\psi_1}$ to $\mathcal{M}^{p,w}_{\psi_2}.$

Next, we assume that $p=1,$ and $f \in \mathcal{M}^{p,w}_{\psi_1}.$ By using Theorem \ref{thm3.2} and the hypothesis concerning $\psi_1$
and $\psi_2,$ we get
$$
\begin{aligned}
\|Mf\|_{W\mathcal{M}^{1,w}_{\psi_2}}
& = \sup_{a\in \mathbb{R}^n, r>0} \frac{1}{\psi_2(a,r)} w(B(a,r))^{-1} \|Mf\|_{WL^{1,w}(B(a,r))} \\
&\leq C \sup_{a\in \mathbb{R}^n, r>0} \frac{1}{\psi_2(a,r)} \sup_{r < t < \infty} w(B(a,t))^{-1} \|f\|_{L^{1,w} B(a,t)}\\
&= C \sup_{a\in \mathbb{R}^n, r>0} \frac{1}{\psi_2(a,r)} \sup_{r < t <\infty} \frac{\psi_1(a,t)}{\psi_1(a,t)} w(B(a,t))^{-1}
\|f\|_{L^{1,w}(B(a,t))} \\
& \leq C \sup_{a\in \mathbb{R}^n, r>0} \frac{1}{\psi_2(a,r)} \sup_{r < t < \infty} \psi_1 (a,t) \|f\|_{\mathcal{M}^{1,w}_{\psi_1}} \\
&= C \|f\|_{\mathcal{M}^{1,w}{\psi_1}} \sup_{a \in \mathbb{R}^n, r>0} \frac{1}{\psi_2(a,r)} \sup_{r < t < \infty} \psi_1(a,t) \\
&\leq C \|f\|_{\mathcal{M}^{1,w}_{\psi_1}}.
\end{aligned}
$$
Therefore, we conclude that $M$ is bounded from $\mathcal{M}^{1,w}_{\psi_1}$ to $W\mathcal{M}^{1,w}_{\psi_1},$
and this completes the proof of Theorem \ref{thm3.3}.\qed

\medskip

Let $a\in\mathbb{R}^n$. Consider the function $t \mapsto \psi_1(a,t)$ on $(0,\infty)$ by $\psi(a,t) = n$ where $t = 1/n$
for some $n\in\mathbb{N}$ and  $\psi(a,t) = te^{-t}$ for otherwise. We also consider the function $t\mapsto\psi_2(a,t)$
by $\psi_2(a,t) = e^{-t}$ for every $t>0$. We can see that
$$
\int_r^\infty \psi_1(a,t) \frac{dt}{t} \leq \psi_2(a,r)
$$
for every $(a,r) \in \mathbb{R}^n \times (0,\infty)$, but there is no $C>0$ such that
$$
\sup_{r<t<\infty} \psi_1(a,t) \leq \psi_2 (a,r)
$$
for every $(a,r)\in\mathbb{R}^n\times (0,\infty).$ Hence, we investigate the condition
$$
\int_{r}^{\infty} \psi_1 (a,t) \frac{dt}{t} \leq C \psi_2 (a,r)
$$
for the boundedness of the Hardy-Littlewood maximal operator and obtain the following results.

\begin{theorem}\label{thm3.4}
Let $1 \leq p < \infty$ and $w \in A_p$. Suppose that $\psi_1$ and $\psi_2$ are nonnegative functions on $\mathbb{R}^n \times (0,\infty)$
satisfying
$$
\int_r^\infty \psi_1 (a,t) \frac{dt}{t} \leq C \psi_2 (a,r)
$$
for every $(a,r)\in \mathbb{R}^n \times (0, \infty)$ where $C$ is a constant that does not depend on $a$ and $r.$ Then,
\begin{enumerate}[(1)]
	\item $M$ is bounded from $\mathcal{M}^{p,w}_{\psi_1}$ to $\mathcal{M}^{p,w}_{\psi_2}$ for $1<p<\infty.$
	\item $M$ is bounded from $\mathcal{M}^{1,w}_{\psi_1}$ to $W\mathcal{M}^{1,w}_{\psi_2}.$
\end{enumerate}
\end{theorem}

Before we present the proof of Theorem \ref{thm3.4}, we prove the following theorem.

\begin{theorem}\label{thm3.5}
Let $1\leq p< \infty$ and $w\in A_p$. Then, for every $a\in\mathbb{R}^n$ and $r>0$,
$$
\|Mf\|_{L^{p,w}(B(a,r))} \leq C_1 w(B(a,r))^\frac1p \int_{r}^\infty w(B(a,s))^{-\frac1p} \|f\|_{L^{p,w}(B(a,s))} \frac{ds}{s},
\quad f\in L_{loc}^{p,w},
$$
for $1<p<\infty$, and
$$
\|Mf\|_{WL^{1,w}(B(a,r))} \leq C_2 w(B(a,r)) \int_{r}^\infty w(B(a,s))^{-1} \|f\|_{L^{1,w}(B(a,s))} \frac{ds}{s}, \quad f\in L_{loc}^{1,w},
$$
where $C_1$ and $C_2$ are positive constants that are independent of $f,$ $a$, and $r.$
\end{theorem}

\textit{Proof.} Given $a\in \mathbb{R}^n$ and $r>0$, we write $f:=f_1+f_2$ where $f_1:=f\cdot \mathcal{X}_{B(a,2r)}$.
Then, by Theorem \ref{garcia}, Theorem \ref{thm3.1}, and Lemma \ref{lem2.8},
$$
\|Mf_1\|_{L^{p,w}(B(a,r))} \leq C \|f\|_{L^{p,w}(B(a,r))} \leq C w(B(a,r))^\frac1p \int_{r}^\infty w(B(a,s))^{-\frac1p}
\|f\|_{L^{p,w}(B(a,s))} \frac{ds}{s}
$$
for $1<p<\infty.$ Meanwhile, for $p=1$, we have
$$
\|Mf_1\|_{WL^{1,w}(B(a,r))} \leq C \|f\|_{L^{1,w}(B(a,r))} \leq C w(B(a,r)) \int_{r}^\infty w(B(a,s))^{-1} \|f\|_{L^{1,w}(B(a,s))}.
\frac{ds}{s}
$$
Since, for every $x\in B(a,r)$,
$$
Mf_2(x) \leq C \sup_{t>r} \frac{1}{|B(a,2t)|} \int_{B(a,2t)} |f(y)|dy,
$$
we have
$$
Mf_2(x) \leq C \sup_{t>2r} \frac{1}{|B(a,t)|} \int_{B(a,t)} |f(y)| dy.
$$
Hence,
$$
Mf_2(x) \leq C \sup_{t>r} \frac{1}{w(B(a,t))^\frac{1}{p}} \|f\|_{L^{p,w}(B(a,t))},
$$
for $1\leq p <\infty.$ Therefore,
$$
Mf_2(x) \leq C \sup_{t>r} \int_t^\infty w(B(a,s))^{-\frac1p} \|f\|_{L^{p,w}(B(a,s))} \frac{ds}{s} \le
\int_r^\infty w(B(a,s))^{-\frac1p} \|f\|_{L^{p,w}(B(a,s))} \frac{ds}{s},
$$
for every $x \in B(a,r)$. It thus follows that
$$
\|Mf_2\|_{L^{p,w}(B(a,r))} \leq C_1 w(B(a,r))^\frac1p \int_r^\infty w(B(a,s))^{-\frac1p} \|f\|_{L^{p,w}(B(a,s))} \frac{ds}{s}
$$
for $1<p<\infty$ and
$$
\|Mf_2\|_{WL^{1,w}(B(a,r))} \leq C_2 w(B(a,r)) \int_r^\infty w(B(a,s))^{-1} \|f\|_{L^{1,w}(B(a,s))} \frac{ds}{s}.
$$
This proves Theorem \ref{thm3.5}.
\qed

\medskip

\textit{Proof of Theorem \ref{thm3.4}}. First, we assume that $1<p<\infty.$ Given $f\in \mathcal{M}_{\psi_1}^{p,w}.$ By using Theorem
\ref{thm3.5} and the assumption concerning $\psi_1$ and $\psi_2$, we get
$$
\begin{aligned}
	\|Mf\|_{\mathcal{M}_{\psi_2}^{p,w}} & = \sup_{a\in \mathbb{R}^n, r>0} \frac{1}{\psi_2(a,r)} \left(\frac{1}{w(B(a,r))}
\int_{B(a,r)} |Mf(x)|^p w(x) dx\right)^{\frac{1}{p}} \\
	& = \sup_{a\in \mathbb{R}^n, r>0} \frac{1}{\psi_2(a,r)} w(B(a,r))^{-\frac{1}{p}} \|Mf\|_{L^{p,w}(B(a,r))} \\
	&\leq C \sup_{a\in \mathbb{R}^n, r>0} \frac{1}{\psi_2(a,r)} \int_r^\infty w(B(a,s))^{-\frac{1}{p}} \|f\|_{L^{p,w}(B(a,s))} \frac{ds}{s}\\
	& = C \sup_{a\in \mathbb{R}^n, r>0} \frac{1}{\psi_2(a,r)} \int_r^\infty \frac{\psi_1(a,s)}{\psi_1(a,s)} w(B(a,s))^{-\frac{1}{p}}
\|f\|_{L^{p,w}(B(a,s))} \frac{ds}{s} \\
	& \leq C \sup_{a \in \mathbb{R}^n, r>0} \frac{1}{\psi_2(a,r)} \int_r^\infty \psi_1 (a,s) \|f\|_{\mathcal{M}_{\psi_1}^{p,w}} \frac{ds}{s} \\
	&= C \|f\|_{\mathcal{M}_{\psi_1}^{p,w}} \sup_{a \in \mathbb{R}^n, r>0} \frac{1}{\psi_2(a,r)} \int_r^\infty \psi_1(a,s) \frac{ds}{s}\\
	&\leq C \|f\|_{\mathcal{M}_{\psi_1}^{p,w}}.
\end{aligned}
$$
Therefore, we conclude that $M$ is bounded from $\mathcal{M}_{p,w}^{\psi_1}$ to $\mathcal{M}_{p,w}^{\psi_2}.$ Next, we assume that $p=1.$
Given $f \in \mathcal{M}_{\psi_1}^{1,w}.$ By using Theorem 3.4 and the assumption concerning $\psi_1$ and $\psi_2,$ we get
$$
\begin{aligned}
	\|Mf\|_{W\mathcal{M}_{\psi_2}^{1,w}} & = \sup_{a\in \mathbb{R}^n, r>0} \frac{1}{\psi_2(a,r)} w(B(a,r))^{-1} \|Mf\|_{WL^{1,w}(B(a,r))} \\
	&\leq C \sup_{a\in \mathbb{R}^n, r>0} \frac{1}{\psi_2(a,r)} \int_r^\infty w(B(a,s))^{-1} \|f\|_{L^{1,w}(B(a,s))} \frac{ds}{s}\\
	&= C \sup_{a\in \mathbb{R}^n, r>0} \frac{1}{\psi_2(a,r)} \int_r^\infty \frac{\psi_1(a,s)}{\psi_1(a,s)} w(B(a,s))^{-1} \|f\|_{L^{1,w}(B(a,s))}
\frac{ds}{s} \\
	& \leq C \sup_{a \in \mathbb{R}^n} \frac{1}{\psi_2(a,r)} \int_r^\infty \psi_1 (a,s) \|f\|_{\mathcal{M}_{\psi_1}^{1,w}} \frac{ds}{s} \\
	&= C \|f\|_{\mathcal{M}_{\psi_1}^{1,w}} \sup_{a \in \mathbb{R}^n, r>0} \frac{1}{\psi_2(a,r)} \int_r^\infty \psi_1(a,s) \frac{ds}{s} \\
	&\leq C \|f\|_{\mathcal{M}_{\psi_1}^{1,w}}.
\end{aligned}
$$
Therefore, we conclude that $M$ is bounded from $\mathcal{M}_{\psi_1}^{1,w}$ to $W\mathcal{M}_{\psi_2}^{1,w}.$ This completes the proof of
Theorem \ref{thm3.4}.
\qed

\section{Fractional Integral and Fractional Maximal Operators on Generalized Weighted Morrey Spaces}

In this section, we prove the boundedness of the fractional integral operator $I_\alpha$ on generalized weighted Morrey spaces
and generalized weighted weak Morrey spaces. The results then imply the boundedness of the fractional maximal operators on
those spaces. As an extention of Theorem \ref{thm4.1} \cite{Wheeden}, we have the following result.

\begin{theorem}\label{thm4.2}
Let $0<\alpha<n, 1\le p<n/\alpha, 1/q = 1/p - \alpha/n$. Then, for every $a\in \mathbb{R}^n$ and $r>0$
$$
\|I_\alpha f\|_{L^{p,w^p}(B(a,r))} \leq C_1 w^q(B(a,r))^{\frac{1}{q}} \int_r^\infty w^q(B(a,s))^{-\frac{1}{q}} \|f\|_{L^{p,w^p}(B(a,s))}
\frac{ds}{s},\quad f\in L_{loc}^{p,w^p}
$$
where $1<p<\infty$, and
$$
\|I_\alpha f\|_{WL^{1,w}(B(a,r))} \leq C_2 w^q(B(a,r))^{\frac{1}{q}} \int_r^\infty w^q(B(a,s))^{-\frac{1}{q}} \|f\|_{L^{1,w}(B(a,s))}
\frac{ds}{s},\quad f\in L_{loc}^{1,w}
$$
where $C_1$ and $C_2$ are constants that do not depend of $f, a,$ and $r.$
\end{theorem}

\textit{Proof.} Given $a\in \mathbb{R}^n$ and $r>0$, we decompose the function $f$ as $f:=f_1+f_2$ where
$f_1 := f\mathcal{X}_{B(a,2r)},$ so that
$$
I_\alpha f (x) = I_\alpha f_1 (x) + I_\alpha f_2 (x).
$$
First, we assume that $1<p<n/\alpha.$ By Theorem \ref{thm4.1}, $I_\alpha$ is bounded from $L^{p,w^p}$ to $L^{q,w^q}.$
Hence,
$$
\|I_\alpha f_1\|_{L^{q,w^q}(B(a,r))} \leq \|I_\alpha f_1\|_{L^{q,w^q}} \leq  C \|f_1\|_{L^{p,w^p}} = C \|f\|_{L^{p,w^p}(B(a,2r))}.
$$
Since $w\in A_{p,q},$ it follows from Theorem \ref{thm2.5} that $w^q \in A_q.$ We see that the map
$r \mapsto \|f\|_{L^{p,w^p}(B(a,2r))}$ is increasing for each $a\in\mathbb{R}^n$, and so by Theorem \ref{garcia} and Lemma \ref{lem2.8} we have
$$
\|I_\alpha f_1\|_{L^{q,w^q}} \leq C w^q(B(a,r))^{\frac{1}{q}} \int_r^\infty w^q(B(a,s))^{-\frac{1}{q}} \|f\|_{L^{p,w^p}(B(a,s))} \frac{ds}{s}.
$$
Next, we obtain the same estimate for $I_\alpha f_2.$ For this, we observe that
$$
|I_\alpha f_2 (x)| \leq \int_{B(a,2r)^c} \frac{|f(y)|}{|x-y|^{n-\alpha}} dy.
$$
The inequalities $|a-x|<r$ and $|x-y| \geq 2r$ implies
$$
\frac{1}{2} |a-y| \leq |x-y| \leq \frac{3}{2} |a-y|.
$$
Then,
$$
|I_\alpha f_2 (x)| \leq C \int_{B(a,2r)^c} \frac{|f(y)|}{|a-y|^{n-\alpha}} dy, \quad x \in B(a,r).
$$
By Fubini's theorem,
$$
\begin{aligned}
 |I_\alpha f_2 (x)|  & \leq C \int_{B(a,2r)^c} \frac{|f(y)|}{|a-y|^{n-\alpha}} dy \\
 & = C \int_{B(a,2r)^c} |f(y)| \int_{|a-y|}^\infty \frac{1}{s^{n-\alpha}} \frac{ds}{s} dy \\
 & = C\int_{r}^\infty \int_{B(a,s) \setminus B(a,r)} \frac{1}{|B(a,s)|^{1-\frac{\alpha}{n}}} |f(y)| dy \frac{ds}{s} \\
 & = \int_r^\infty \frac{1}{|B(a,s)|^{1+\frac{1}{q}-\frac{1}{p}}} \int_{B(a,s)} |f(y)| dy \frac{ds}{s}.
\end{aligned}
$$
Next, by H\"older's inequality and the assumption that $w\in A_{p,q}$, we have
$$
\begin{aligned}
	&\frac{1}{|B(a,s)|^{1+\frac{1}{q}-\frac{1}{p}}} \int_{B(a,s)} |f(y)| dy \\
& = \frac{1}{|B(a,s)|^{1 + \frac{1}{q}-\frac{1}{p}}} \int_{B(a,s)} \frac{|f(y)|w(y)}{w(y)} dy \\
& \leq \frac{1}{|B(a,s)|^{1+\frac{1}{q}-\frac{1}{p}}} \left(\int_{B(a,s)} |f(y)|^p w(y)^p dy \right)^{\frac{1}{p}}
\left( \int_{B(a,s)} w(y)^{-p'} dy \right)^{\frac{1}{p'}} \\
& = w^q(B(a,s))^{-\frac{1}{q}}\|f\|_{L^{p,w^p}(B(a,s))} \left(\frac{1}{|B(a,s)|} \int_{B(a,s)} w(y)^q dy\right)^{\frac{1}{q}}
\left(\frac{1}{|B(a,s)|} \int_{B(a,s)} w(y)^{-p'} dy\right)^{\frac{1}{p'}} \\
& \leq C w^q(B(a,s))^{-\frac{1}{q}}\|f\|_{L^{p,w^p}(B(a,s))}.
\end{aligned}
$$
Hence,
$$
|I_\alpha f_2 (x)| \leq C \int_r^\infty w^q(B(a,s))^{-\frac{1}{q}}\|f\|_{L^{p,w^p}(B(a,s))} \frac{ds}{s}, \quad x \in B(a,r),
$$
and this implies that
$$
\|I_\alpha f_2\|_{L^{q,w^q}(B(a,r))} \leq C w^q(B(a,r))^{\frac{1}{q}} \int_r^\infty w^q(B(a,s))^{-\frac{1}{q}}\|f\|_{L^{p,w^p}(B(a,s))}
\frac{ds}{s}.
$$
Therefore,
$$
\|I_\alpha f\|_{L^{q,w^q}(B(a,r))} \leq C_1 w^q(B(a,r))^{\frac{1}{q}} \int_r^\infty w^q(B(a,s))^{-\frac{1}{q}}\|f\|_{L^{p,w^p}(B(a,s))}
\frac{ds}{s}.
$$

Next, we assume $p=1.$ Then,
$$
\|I_\alpha f\|_{WL^{1,w}(B(a,r))} \leq\|I_\alpha f_1\|_{WL^{1,w}(B(a,r))}+\|I_\alpha f_2\|_{WL^{1,w}(B(a,r))}.
$$
By Theorem \ref{thm4.1}, we have
$$
\|I_\alpha f_1\|_{WL^{1,w}(B(a,r))} \leq \|I_\alpha f_1\|_{WL^{1,w}} \leq C \|I_\alpha f_1\|_{L^{1,w}} =
C \|I_\alpha f\|_{L^{1,w}(B(a,2r))}.
$$
Since $w\in A_{1,q},$ it follows from Theorem \ref{thm2.5} that $w^q\in A_q.$ By the same argument as for the case $p>1,$
we have that
$$
\|I_\alpha f\|_{L^{1,w}(B(a,r))} \leq C w^q(B(a,r))^{\frac{1}{q}} \int_r^\infty w^q(B(a,s))^{-\frac{1}{q}} \|f\|_{L^{1,w}(B(a,s))}
\frac{ds}{s}.
$$
Since
$$
|I_\alpha f_2 (x)| \leq C \int_{B(a,2r)^c} \frac{|f(y)|}{|a-y|^{n-\alpha}} dy, \quad x \in B(a,r),
$$
by Fubini's theorem we have
$$
\begin{aligned}
|I_\alpha f_2 (x)|
& \leq C \int_{\mathbb{R}^n \setminus B(a,2r)} \frac{|f(y)|}{|a-y|^{n-\alpha}} dy \\
& = C \int_{B(a,2r)^c} |f(y)| \int_{|a-y|}^\infty \frac{1}{s^{n-\alpha}} \frac{ds}{s} dy \\
& = C \int_{r}^\infty \int_{B(a,s) \setminus B(a,r)} \frac{1}{|B(a,s)|^{1-\frac{\alpha}{n}}} |f(y)| dy \frac{ds}{s} \\
& = \int_r^\infty \frac{1}{|B(a,s)|^{\frac{1}{q}}} \int_{B(a,s)} |f(y)| dy \frac{ds}{s}.
\end{aligned}
$$
By H\"older's inequality, together with the assumption that $w\in A_{1,q}$ and the fact that $q > 1$, we have
$$
\begin{aligned}
& \frac{1}{|B(a,s)|^{\frac{1}{q}}} \int_{B(a,s)} |f(y)| dy\\
& \leq \frac{1}{|B(a,s)|^{\frac{1}{q}}} \int_{B(a,s)} |f(y)| \frac{w(y)}{w(y)} dy \\
& \leq \frac{1}{|B(a,s)|^{\frac{1}{q}}} \left(\int_{B(a,s)} |f(y)|w(y) dy\right) \|w^{-1}\|_{L^\infty (B(a,s))} \\
& = w^q(B(a,s))^{-\frac{1}{q}} \|f\|_{L^{1,w}(B(a,s))} \left(\frac{1}{|B(a,s)|} \int_{B(a,s)} w(x)^q \right)^{\frac{1}{q}}
\|w^{-1}\|_{L^\infty (B(a,s))} \\
& \leq w^q(B(a,s))^{-\frac{1}{q}} \|f\|_{L^{1,w}(B(a,s))} \left(\frac{1}{|B(a,s)|} \int_{B(a,s)} w(x)^q \right)
\|w^{-1}\|_{L^\infty (B(a,s))} \\
& \leq C  w^q(B(a,s))^{-\frac{1}{q}} \|f\|_{L^{1,w}(B(a,s))}.
\end{aligned}
$$
Hence, for $x \in B(a,r)$,
$$
\begin{aligned}
|I_\alpha f_2(x)|
& \leq C \int_r^\infty \frac{1}{|B(a,s)|^{\frac{1}{q}}} \int_{B(a,s)} |f(y)| dy \frac{ds}{s} \\
& \leq C \int_r^\infty w^q(B(a,s))^{-\frac{1}{q}} \|f\|_{L^{1,w}(B(a,s))} \frac{ds}{s},
\end{aligned} $$
and H\"older's inequality implies that
$$
\begin{aligned}
& \|I_\alpha f_2\|_{WL^{1,w}(B(a,r))} \\
& = \sup_{\gamma > 0} \gamma w(\{x\in B(a,r):|I_\alpha f_2 (x)|>\gamma\}) \\
&= \sup_{\gamma>0} \gamma \int_{\{x\in B(a,r):|I_\alpha f_2 (x)|>\gamma\}} w(x) dx \\
& \leq \sup_{\gamma>0} \gamma \left(\int_{\{x\in B(a,r):|I_\alpha f_2 (x)|>\gamma\}} |I_\alpha f_2 (y)| w(y) dy\right)
\left\|\frac{1}{I_\alpha f_2}\right\|_{L^\infty (\{x\in B(a,r):|I_\alpha f_2 (x)|>\gamma\})} \\
&= \sup_{\gamma>0} \gamma \left(\int_{B(a,r)} |I_\alpha f_2 (y)| w(y) dy\right) \frac{1}{\gamma} \\
&= \sup_{\gamma>0} \left(\int_{B(a,r)} |I_\alpha f_2 (y)| w(y) dy\right) \\
& \leq C \sup_{\gamma>0} \int_{B(a,r)} \left[\int_r^\infty w^q(B(a,s))^{-\frac{1}{q}}
\|f\|_{L^{1,w}(B(a,s))} \frac{ds}{s} \right] w(y) dy \\
& = C \left(\int_{B(a,r)} w(y) dy \right) \left[\int_r^\infty w^q(B(a,s))^{-\frac{1}{q}} \|f\|_{L^{1,w}(B(a,s))} \frac{ds}{s} \right] \\
& \leq C \left(\int_{B(a,r)} w(y)^q dy \right)^{\frac{1}{q}} \left[\int_r^\infty w^q(B(a,s))^{-\frac{1}{q}} \|f\|_{L^{1,w}(B(a,s))}
\frac{ds}{s} \right] \\
& = C w^q(B(a,r))^{\frac{1}{q}} \int_r^\infty w^q(B(a,s))^{-\frac{1}{q}} \|f\|_{L^{1,w}(B(a,s))} \frac{ds}{s}.
\end{aligned}
$$
Therefore,
$$
\|I_\alpha f\|_{WL^{1,w}(B(a,r))} \leq C_2 w^q(B(a,r))^{\frac{1}{q}} \int_r^\infty w^q(B(a,s))^{-\frac{1}{q}} \|f\|_{L^{1,w}(B(a,s))}
\frac{ds}{s}.
$$
This completes the proof of Theorem \ref{thm4.2}. \qed

\begin{theorem}\label{thm4.3}
Let $0<\alpha<n, 1 \leq p < n/\alpha, 1/q = 1/p -\alpha/n, w \in A_{p,q}$, and $I_\alpha$ be the fractional integral operator.
Suppose that $\psi_1$ and $\psi_2$ are nonnegative functions on $\mathbb{R}^n \times (0,\infty)$ satisfying
$$
\int_r^\infty \frac{w^p(B(a,t))^{\frac{1}{p}}}{w^q(B(a,t))^{\frac{1}{q}}} \psi_1 (a,t) \frac{dt}{t} \leq C \psi_2 (a,r)
$$
for every $(a,r)\in \mathbb{R}^n \times (0, \infty)$, where $C$ is a constant that does not depend on $a$ and $r.$ Then,
\begin{enumerate}[(1)]
	\item $I_\alpha$ is bounded from $\mathcal{M}^{p,w^p}_{\psi_1}$ to $\mathcal{M}^{q,w^q}_{\psi_2}$ for $1<p<\infty.$
	\item $I_\alpha$ is bounded from $\mathcal{M}^{1,w}_{\psi_1}$ to $W\mathcal{M}^{q,w^q}_{\psi_2}.$
\end{enumerate}
\end{theorem}

\textit{Proof} First, we assume that $1<p<n/\alpha$, and let $f\in \mathcal{M}^{p,w}_{\psi_1}.$ By using Theorem \ref{thm4.2} and
the assumption on $\psi_1$ and $\psi_2,$ we get
$$
\begin{aligned}
\|I_\alpha f\|_{\mathcal{M}^{q,w^q}_{\psi_2}}
&= \sup_{a\in \mathbb{R}^n, r>0} \frac{1}{\psi_2(a,r)} \left(\frac{1}{w^q(B(a,r))} \int_{B(a,r)}
|I_\alpha f(x)|^q w(x)^q dx\right)^{\frac{1}{q}} \\
&= \sup_{a\in \mathbb{R}^n, r>0} \frac{1}{\psi_2(a,r)} w^q(B(a,r))^{-\frac{1}{q}} \|I_\alpha\|_{L^{p,w^p}(B(a,r))} \\
&\leq C \sup_{a\in \mathbb{R}^n, r>0} \frac{1}{\psi_2(a,r)} \int_r^\infty w^q(B(a,s))^{-\frac{1}{q}} \|f\|_{L^{p,w^p}(B(a,s))} \frac{ds}{s}\\
&= C \sup_{a\in \mathbb{R}^n, r>0} \frac{1}{\psi_2(a,r)} \int_r^\infty \frac{\psi_1(a,s)}{\psi_1(a,s)} w^q(B(a,s))^{-\frac{1}{q}}
\|f\|_{L^{p,w^p}(B(a,s))} \frac{ds}{s} \\
&\leq C \frac{1}{\psi_2(a,r)} \int_r^\infty \psi_1 (a,s) \frac{w^q(B(a,s))^{-\frac{1}{q}}}{w^p(B(a,s))^{-\frac{1}{p}}}
\|f\|_{\mathcal{M}^{p,w^p}_{\psi_1}} \frac{ds}{s} \\
&= C \|f\|_{\mathcal{M}^{p,w^p}_{\psi_1}} \sup_{a \in \mathbb{R}^n, r>0} \frac{1}{\psi_2(a,r)} \int_t^\infty \psi_1(a,s)
\frac{w^p(B(a,s))^{\frac{1}{p}}}{w^q(B(a,s))^{\frac{1}{q}}} \frac{ds}{s} \\
&\leq C \|f\|_{\mathcal{M}^{p,w^p}_{\psi_1}}.
\end{aligned}
$$
Therefore, we conclude that $I_\alpha$ is bounded from $\mathcal{M}^{p,w}_{\psi_1}$ to $\mathcal{M}^{p,w}_{\psi_2}.$

Next, we assume that $p=1,$ and let $f \in \mathcal{M}^{p,w}_{\psi_1}.$ By using Theorem \ref{thm4.2} and the assumption on
$\psi_1$ and $\psi_2,$ we get
$$
\begin{aligned}
\|I_\alpha f\|_{W\mathcal{M}^{q,w^q}_{\psi_2}}
& = \sup_{a\in \mathbb{R}^n, r>0} \frac{1}{\psi_2(a,r)} w^q(B(a,r))^{-\frac{1}{q}} \|I_\alpha\|_{WL^{1,w}(B(a,r))} \\
&\leq C \sup_{a\in \mathbb{R}^n, r>0} \frac{1}{\psi_2(a,r)} \int_r^\infty w^q(B(a,s))^{-\frac{1}{q}} \|f\|_{L^{p,w^p}(B(a,s))}
\frac{ds}{s}\\
&= C \sup_{a\in \mathbb{R}^n, r>0} \frac{1}{\psi_2(a,r)} \int_r^\infty \frac{\psi_1(a,s)}{\psi_1(a,s)} w^q(B(a,s))^{-\frac{1}{q}}
\|f\|_{L^{1,w}(B(a,s))} \frac{ds}{s} \\
& \leq C \frac{1}{\psi_2(a,r)} \int_r^\infty \psi_1 (a,s) \frac{w^q(B(a,s))^{-\frac{1}{q}}}
{w(B(a,s))^{-1}} \|f\|_{\mathcal{M}^{1,w}_{\psi_1}} \frac{ds}{s} \\
&= C \|f\|_{\mathcal{M}^{1,w}_{\psi_1}} \sup_{a \in \mathbb{R}^n, r>0} \frac{1}{\psi_2(a,r)} \int_t^\infty \psi_1(a,s)
\frac{w(B(a,s))}{w^q(B(a,s))^{\frac{1}{q}}} \frac{ds}{s} \\
&\leq C \|f\|_{\mathcal{M}^{1,w}_{\psi_1}}.
\end{aligned}
$$
Therefore, we conclude that $I_\alpha$ is bounded from $\mathcal{M}^{1,w}_{\psi_1}$ to $W\mathcal{M}^{1,w}_{\psi_1}.$
This completes the proof of Theorem \ref{thm4.3}. \qed

\medskip

By the relation (\ref{Eq1}) and Theorem \ref{thm4.3}, we have the following corollary for the fractional maximal operator $M_\alpha$.

\begin{cor}
Let $0<\alpha<n, 1 \leq p < n/\alpha, 1/q = 1/p -\alpha/n$, and $w \in A_{p,q}$. Suppose that $\psi_1$ and $\psi_2$ are
nonnegative functions on $\mathbb{R}^n \times (0,\infty)$ satisfying
$$
\int_r^\infty \frac{w^p(B(a,t))^{\frac{1}{p}}}{w^q(B(a,s))^{\frac{1}{q}}} \psi_1(a,s) \frac{dt}{t} \leq C \psi_2(a,r)
$$
for every $(a,r)\in \mathbb{R}^n \times (0, \infty)$ where $C$ is a constant that does not depend on $r.$ Then,
\begin{enumerate}[(1)]
	\item $M_\alpha$ is bounded from $\mathcal{M}^{p,w^p}_{\psi_1}$ to $\mathcal{M}^{q,w^q}_{\psi_2}$ for $1<p<\infty.$
	\item $M_\alpha$ is bounded from $\mathcal{M}^{1,w}_{\psi_1}$ to $W\mathcal{M}^{q,w^q}_{\psi_2}.$
\end{enumerate}
\end{cor}

\section{Calderon-Zygmund Operators on Weighted Morrey Spaces}

In this section, we prove the boundedness of the Calderon-Zygmund operators $T=T_K$ on generalized weighted Morrey spaces
and generalized weighted weak Morrey spaces. As stated earlier, we have Theorem \ref{thm5.1} about the boundedness
of the Calderon-Zygmund operators on weighted Lebesgue spaces and weighted weak Lebesgue spaces. As an extention of this theorem,
we have the following result.

\begin{theorem}\label{thm5.2}
Let $1 \leq p < \infty, w\in A_p$ and $f\in L_{loc}^{p,w}.$ Then, for every $a\in \mathbb{R}^n$ and $r>0$
$$
\|Tf\|_{L^{p,w} (B(a,r))} \leq C_1 w(B(a,r))^{\frac{1}{p}} \int_r^\infty w(B(a,s))^{-\frac{1}{p}} \|f\|_{L^{p,w}(B(a,s))} \frac{ds}{s}
$$
where $1<p<\infty$, and
$$
\|Tf\|_{WL^{1,w} (B(a,r))} \leq C_1 w(B(a,r))\int_r^\infty w(B(a,s))^{-1} \|f\|_{L^{p,w}(B(a,s))} \frac{ds}{s},
$$
where $C_1$ and $C_2$ are constants that do not depend on $f, a,$ and $r.$
\end{theorem}

\textit{Proof.} Write $f$ as $f:=f_1+f_2$ where $f_1:=f\cdot \mathcal{X}_{B(a,2r)}.$ First, we consider the case where $1<p<\infty.$
Then, since $w\in A_p$, we know that $T$ is bounded on $L^{p,w}.$ Thus, for every $a\in \mathbb{R}^n$ and $r>0$, we have
$$
\|Tf\|_{L^{p,w}(B(a,r))} \leq \|Tf_1\|_{L^{p,w}(B((a,r))} + \|Tf_1\|_{L^{p,w}(B((a,r))}
$$
and
$$
\|Tf_1\|_{L^{p,w}(B(a,r))} \leq \|Tf_1\|_{L^{p,w}} \leq C \|f_1\|_{L^{p,w}} \leq C \|f\|_{L^{p,w}(B(a,2r))}.
$$
By Theorem \ref{garcia} and Lemma \ref{lem2.8},
$$
\|Tf_1\|_{L^{p,w}(B(a,r))} \leq C w(B(a,r))^{\frac{1}{p}} \int_r^\infty \frac{1}{w(B(a,s))^{\frac{1}{p}}}
\|f\|_{L^{p,w}(B(a,s))} \frac{ds}{s}.
$$
Note that for $x\in B(a,r)$, we have
$$
|Tf_2(x)| \leq C \int_{B(a,2r)^c} \frac{|f(y)|}{|x-y|^n} dy.
$$
On other hand, the inequalities $|a-x|<r$ and $|x-y|\geq 2r$ imply that
$$
\frac{1}{2} |a-y|\leq |x-y| \leq \frac{3}{2} |a-y|.
$$
Then, by using Lemma \ref{lem2.9},
$$
|Tf_2(x)| \leq C \int_{B(a,2r)^c} \frac{|f(y)|}{|a-y|^n} dy \leq C \int_r^\infty w(B(a,s))^{-\frac{1}{p}}
\|f\|_{L^{p,w}(B(a,s))} \frac{ds}{s}.
$$
Hence,
$$
\|Tf_2\|_{L^{p,w}(B(a,r))} \leq C w(B(a,r))^{\frac{1}{p}} \int_r^\infty w(B(a,s))^{-\frac{1}{p}} \|f\|_{L^{p,w}(B(a,s))} \frac{ds}{s}.
$$
Therefore, we conclude that
$$
\|Tf\|_{L^{p,w}(B(a,r))} \leq C_1 w(B(a,r))^{\frac{1}{p}} \int_r^\infty w(B(a,s))^{-\frac{1}{p}}\|f\|_{L^{p,w}(B(a,s))} \frac{ds}{s}.
$$

Next, for the case where $p=1,$ we have
$$
\|Tf\|_{WL^{1,w}(B(a,r))} \leq \|Tf_1\|_{WL^{1,w}(B(a,r))} + \|Tf_1\|_{WL^{1,w}(B(a,r))}
$$
and
$$
\|Tf_1\|_{WL^{1,w}(B(a,r))} \leq \|Tf_1\|_{L^{1,w}} \leq C \|f_1\|_{L^{1,w}} \leq C \|f\|_{L^{1,w}(B(a,2r))}
$$
for every $a\in \mathbb{R}^n$ and $r>0$. By the boundedness of $T$ from $L^{1,w}$ to $WL^{1,w},$ we have
$$
\|Tf_1\|_{WL^{1,w}(B(a,r))} \leq \|Tf_1\|_{WL^{1,w}} \leq C \|f_1\|_{L^{1,w}} \leq C \|f\|_{L^{1,w}(B(a,2r))}.
$$
By Theorem \ref{garcia} and Lemma \ref{lem2.8},
$$
\|Tf_1\|_{WL^{1,w}(B(a,r))} \leq C w(B(a,r)) \int_r^\infty w(B(a,s))^{-1} \|f\|_{L^{p,w}(B(a,s))} \frac{ds}{s}.
$$
As we do before, for $x\in B(a,r)$ we have
$$
|Tf_2(x)| \leq C \int_r^\infty w(B(a,s))^{-1} \|f\|_{L^{p,w}(B(a,s))} \frac{ds}{s}.
$$
Thus, by H\"older's inequality,
$$
\begin{aligned}
&\|Tf\|_{WL^{1,w}(B(a,r))} \\
&= \sup_{\gamma > 0} \gamma w(\{x\in B(a,r):|Tf_2(x)|>\gamma\}) \\
&= \sup_{\gamma>0} \gamma \int_{\{x\in B(a,r):|Tf_2(x)|>\gamma\}} w(x) dx \\
&\leq \sup_{\gamma>0} \gamma\int_{\{x\in B(a,r):|Tf_2(x)|>\gamma\}} |Tf_2(x)| w(x) dx \cdot
\left\|\frac{1}{Tf_2} \right\|_{L^\infty (\{x\in B(a,r):|Tf_2(x)|>\gamma\})} \\
&\leq \sup_{\gamma>0} \gamma \int_{\{x\in B(a,r):|Tf_2(x)|>\gamma\}} |Tf_2 (x)|w(x) dx \cdot \frac{1}{\gamma} \\
&\leq \int_{B(a,r)} |Tf_2(x)| w(x) dx \\&\leq C \int_{B(a,r)} \int_r^\infty w(B(a,s))^{-1} \|f\|_{L^{p,w}(B(a,s))} \frac{ds}{s} w(x) dx \\
&= C w(B(a,r)) \int_r^\infty w(B(a,s))^{-1} \|f\|_{L^{1,w}(B(a,s))} \frac{ds}{s}.
\end{aligned}
$$
Therefore,
$$
\|Tf\|_{L^{1,w}(B(a,r))} \leq C_2 w(B(a,r))\int_r^\infty w(B(a,s))^{-1} \|f\|_{L^{p,w}(B(a,s))} \frac{ds}{s},
$$
and this proves Theorem \ref{thm5.2}. \qed

\begin{theorem} \label{thm5.3}
Let $1 \leq p < \infty, w \in A_p$, and $T$ be the Calderon-Zygmund operator. Suppose that $\psi_1$ and $\psi_2$ are nonnegative functions
on $\mathbb{R}^n \times (0,\infty)$ satisfying
$$
\int_r^\infty \psi_1 (a,t) \frac{dt}{t} \leq C \psi_2 (a,r)
$$
for every $(a,r)\in \mathbb{R}^n \times (0, \infty)$ where $C$ is a constant that does not depend on $a$ and $r.$ Then,
\begin{enumerate}[(1)]
	\item $T$ is bounded from $\mathcal{M}^{p,w}_{\psi_1}$ to $\mathcal{M}^{p,w}_{\psi_2}$ for $1<p<\infty.$
	\item $T$ is bounded from $\mathcal{M}^{1,w}_{\psi_1}$ to $W\mathcal{M}^{1,w}_{\psi_2}.$
\end{enumerate}
\end{theorem}

\textit{Proof} First, we assume that $1<p<\infty.$ Given $f\in \mathcal{M}_{p,w}^{\psi_1}.$ By using Theorem \ref{thm5.2} and the
assumption concerning $\psi_1$ and $\psi_2,$ we get
$$
\begin{aligned}
\|Tf\|_{\mathcal{M}_{\psi_2}^{p,w}}
&= \sup_{a\in \mathbb{R}^n, r>0} \frac{1}{\psi_2(a,r)} \left(\frac{1}{w(B(a,r))}
\int_{B(a,r)} |Tf(x)|^p w(x) dx\right)^{\frac{1}{p}} \\
&= \sup_{a\in \mathbb{R}^n, r>0} \frac{1}{\psi_2(a,r)} w(B(a,r))^{-\frac{1}{p}} \|Tf\|_{L^{p,w}(B(a,r))} \\
&\leq C \sup_{a\in \mathbb{R}^n, r>0} \frac{1}{\psi_2(a,r)} \int_r^\infty w(B(a,s))^{-\frac{1}{p}} \|f\|_{L^{p,w}(B(a,s))} \frac{ds}{s}\\
&= C \sup_{a\in \mathbb{R}^n, r>0} \frac{1}{\psi_2(a,r)} \int_r^\infty \frac{\psi_1(a,s)}{\psi_1(a,s)} w(B(a,s))^{-\frac{1}{p}}
\|f\|_{L^{p,w}(B(a,s))} \frac{ds}{s} \\
&\leq C \sup_{a \in \mathbb{R}^n, r>0} \frac{1}{\psi_2(a,r)} \int_r^\infty \psi_1 (a,s) \|f\|_{\mathcal{M}^{p,w}_{\psi_1}} \frac{ds}{s} \\
&= C \|f\|_{\mathcal{M}^{p,w}_{\psi_1}} \sup_{a \in \mathbb{R}^n, r>0} \frac{1}{\psi_2(a,r)} \int_r^\infty \psi_1(a,s) \frac{ds}{s} \\
&\leq C \|f\|_{\mathcal{M}_{\psi_1}^{p,w}}.
\end{aligned}
$$
Therefore, we conclude that $T$ is bounded from $\mathcal{M}_{p,w}^{\psi_1}$ to $\mathcal{M}^{p,w}_{\psi_2}.$

Next, we assume that $p=1,$ and let $f \in \mathcal{M}^{p,w}_{\psi_1}.$ By using Theorem \ref{thm5.2} and the assumption on $\psi_1$ and
$\psi_2,$ we get
$$
\begin{aligned}
\|Tf\|_{W\mathcal{M}_{\psi_2}^{p,w}}
&= \sup_{a\in \mathbb{R}^n, r>0} \frac{1}{\psi_2(a,r)} w(B(a,r))^{-\frac{1}{p}}
\|Tf\|_{WL^{1,w}(B(a,r))} \\
&\leq C \sup_{a\in \mathbb{R}^n, r>0} \frac{1}{\psi_2(a,r)} \int_r^\infty w(B(a,s))^{-1} \|f\|_{L^{1,w}(B(a,s))} \frac{ds}{s}\\
&= C \sup_{a\in \mathbb{R}^n, r>0} \frac{1}{\psi_2(a,r)} \int_r^\infty \frac{\psi_1(a,s)}{\psi_1(a,s)} w(B(a,s))^{-1}
\|f\|_{L^{1,w}(B(a,s))} \frac{ds}{s} \\
&\leq C \sup_{a \in \mathbb{R}^n} \frac{1}{\psi_2(a,r)} \int_r^\infty \psi_1(a,s)
\|f\|_{\mathcal{M}^{1,w}_{\psi_1}} \frac{ds}{s} \\&= C \|f\|_{\mathcal{M}^{1,w}_{\psi_1}} \sup_{a \in \mathbb{R}^n, r>0} \frac{1}{\psi_2(a,r)}
\int_r^\infty \psi_1(a,s) \frac{ds}{s} \\ &\leq C \|f\|_{\mathcal{M}^{1,w}_{\psi_1}}.
\end{aligned}
$$
Therefore, we conclude that $T$ is bounded from $\mathcal{M}^{1,w}_{\psi_1}$ to $W\mathcal{M}^{1,w}_{\psi_1}.$
This completes the proof of Theorem \ref{thm5.3}. \qed

\begin{remark}
By using the results in \cite{Quek}, we can extend Theorem \ref{thm5.3} by replacing $T$ with $\theta$-type Calderon-Zygmund operators
$T_{\theta}$. The definition of $\theta$-type Calderon-Zygmund Operator $T_{\theta}$ may be found in \cite{Yabuta}. Accordingly one can
obtain a result that is more general than \cite{Wang}.
\end{remark}

\bigskip

\noindent{\bf Acknowledgement}. The second author is supported by P2MI-ITB 2022 Program.

\end{document}